\documentclass[11pt,twoside]{smfart}
\usepackage{a4}
\usepackage[dvips]{graphicx}
\usepackage[francais,english]{babel}
\usepackage{amsmath}
\usepackage{amsthm}
\usepackage{amssymb}
\usepackage{amsfonts}
\usepackage{mathrsfs}
\usepackage[all]{xy}
\theoremstyle{plain}
\newtheorem{theo}[subsection]{Th\'eor\`eme}

\newtheorem{reve}[subsection]{R\^eve}

\theoremstyle{definition}
\newtheorem{defin}[subsection]{D\'efinition}

\theoremstyle{remark}

\newcommand{\Oc}{\mathcal{O}}

\newcommand{\Sc}{\mathcal{S}}

\newcommand{\Rc}{\mathcal{R}}
\newcommand{\Gc}{\mathcal{G}}

\newcommand{\Ib}{\mathbb{I}}
\newcommand{\Db}{\mathbb{D}}
\newcommand{\A}{\mathbb{A}}
\newcommand{\C}{\mathbb{C}}
\newcommand{\R}{\mathbb{R}}
\newcommand{\Q}{{\mathbb{Q}}}

\newcommand{\qb}{{\overline{\Q}}}
\newcommand{\Z}{\mathbb{Z}}
\newcommand{\Nm}{\mathrm{Nm}}

\newcommand{\Tb}{\mathbb{T}}
\newcommand{\Ab}{\mathbb{A}}

\newcommand{\G}{\mathbb{G}}
\newcommand{\Sb}{\mathbb{S}}
\newcommand{\Hb}{\mathbb{H}}
\newcommand{\Pb}{\mathbb{P}}

\newcommand{\nat}{\mathfrak{nat}}
\newcommand{\rec}{\mathfrak{rec}}

\newcommand{\artin}{\mathfrak{artin}}
\newcommand{\Hom}{\mathrm{Hom}}

\newcommand{\im}{\mathrm{im}}
\newcommand{\End}{\mathrm{End}}
\newcommand{\diag}{\mathrm{diag}}

\newcommand{\Ell}{\mathscr{E}\ell\ell}

\newcommand{\limproj}{\underset{\leftarrow}{\mathrm{lim}}\,}
\newcommand{\Gal}{\mathrm{Gal}}

\newcommand{\PGL}{\mathrm{PGL}}

\newcommand{\NR}{\mathrm{NR}}

\newcommand{\HC}{\mathrm{HC}}
\newcommand{\GL}{\mathrm{GL}}

\newcommand{\GSp}{\mathrm{GSp}}

\newcommand{\BMT}{\mathrm{BMT}}

\newcommand{\Sh}{\mathrm{Sh}}

\newcommand{\Res}{\mathrm{Res}}

\title{Quelques bords irrationnels de vari\'et\'es de Shimura}
\author[F. Paugam]{Fr\'ed\'eric Paugam}
\address{NWF I - Mathematik, Universit\"{a}t Regensburg, 93040 Regensburg, 
Germany}
\email{frederic.paugam@mathematik.uni-regensburg.de}
\date{8 juillet 2004}
\makeatletter

\begin{document}

\begin{abstract}
We are looking for a formulation of Manin's real multiplication
question in higher rank. This question comports, in our point of view, at least two steps:
\begin{enumerate}
\item a formalisation of the linear algebra side of the story in terms of morphisms of algebraic
groups analogous to Shimura and Deligne's point of view on the theory of complex multiplication.
\item a work of noncommutative algebraic geometry.
\end{enumerate}
We are interested only by the first step, and recall the know results in litterature on the second
one.
Our point of view, perhaps too naive, is that before looking for a good notion of noncommutative
abelian variety, we have to know what are their periods and what we want to do with them.
The appendix contains additional notes not exposed.
\end{abstract}
\maketitle

%\begin{center}
%\it
%\item Notes \'etay\'ees pour un expos\'e aux seminaires de\\
%\item g\'eom\'etrie alg\'ebrique de Rennes (29/4/04),\\
%\item K-theorie et g\'eom\'etrie non commutative \`a l'IHP, Paris (3/5/04),\\

%\item math\'ematiques pures de l'ENS Lyon (13/5/04),\\
%\item g\'eom\'etrie non commutative et th\'eorie des nombres, II, MPIM-Bonn (6/04),\\
%\item g\'eom\'etrie diophantienne de Goettingen (6/04).
%\end{center}

%\selectlanguage{english}

%\selectlanguage{francais}

%\abstract{On cherche une formulation en rang sup\'erieur de la question de
%multiplication r\'eelle pos\'ee par Manin. Cette question
%comporte, de notre point de vue, deux phases:
%\begin{enumerate}
%\item une formalisation en termes de morphismes de groupes alg\'ebriques
%analogue de la formulation \`a la Shimura/Deligne de la th\'eorie de la
%multiplication complexe.
%\item un travail de g\'eom\'etrie alg\'ebrique non commutative.
%\end{enumerate}
%On s'int\'eresse uniquement \`a la premi\`ere phase, et on rappelle les
%r\'esultats connus dans la litt\'erature sur la seconde.

%On prend le point de vue peut-\^etre trop simpliste que, pour savoir ce
%qu'est une notion utile de vari\'et\'es ab\'eliennes non commutatives,
%il faut au moins conna\^itre leurs p\'eriodes, et ce qu'on peut en faire.
%On met en appendice des notes additionnelles non expos\'ees.}

%************************************************************************
\section{Introduction}
Le th\'eor\`eme de Kronecker-Weber nous dit que l'extension ab\'elienne
maximale de $\Q$ est engendr\'ee par l'image de l'application
$$e^{2i\pi\bullet}:\Q\to\C.$$

Soit $F/\Q$ un corps de nombres.
Le douzi\`eme probl\`eme de Hilbert est de trouver des fonctions analytiques
rempla\c{c}ant la fonction exponentielle et dont les valeurs sp\'eciales
engendrent l'extension ab\'elienne maximale de $F$. Ceci peut-\^etre
consid\'er\'e comme une explicitation de la th\'eorie du corps de classe
qui permet, elle, de d\'ecrire cette extension en termes de classes d'id\'eaux
g\'en\'eralis\'es.

Soit $F/\Q$ un corps quadratique imaginaire. Alors $\C/\Oc_F$ est une courbe
elliptique dite \`a multiplication complexe par $F$ et on peut montrer qu'elle
est d\'efinie sur un corps de nombres, extension ab\'elienne de $F$. La
th\'eorie de la multiplication complexe (``Jugendtraum'' de Kronecker) dit que
toutes les extension ab\'eliennes de $F$ peuvent \^etre d\'ecrites gr\^ace \`a
des courbes elliptiques de ce genre et \`a leurs points de torsion.

Shimura et Taniyama ont montr\'e que le m\^eme genre de chose marchait pour
des extensions quadratiques totalement imaginaires de corps totalement r\'eels
en utilisant des vari\'et\'es ab\'eliennes.

Ce que Manin appelle son ``Alterstraum'' (voir \cite{Manin3}) est que si $F/\Q$ est quadratique
r\'eel, on peut utiliser des ``quotients'' du type $\C/\Oc_F$ pour faire
une th\'eorie de la multiplication r\'eelle analogue \`a ce qui se fait
avec les courbes elliptiques.
On remarque que dans ce cas, $\Oc_F$ est dense dans une
droite r\'eelle et la notion habituelle de quotient ne suffit pas.
Il faut remarquer que Gauss, Shimura, Stark et bien d'autres ont aussi \'etudi\'e
ce genre de r\^eve, sans le formuler sous cette forme.

J'ai cherch\'e \`a mieux comprendre ce r\^eve de Manin, en prenant un point
de vue \`a la Deligne sur les vari\'et\'es de Shimura qui permet plus
facilement de voir ce dont on a besoin en dimension sup\'erieure.

Regardons maintenant le point de vue modulaire sur la multiplication
complexe.
L'espace des modules des courbes elliptiques sur $\C$ \`a structure de niveau
infinie est une vari\'et\'e de Shimura qui \`a un mod\`ele
sur $\Q$ (l'espace de module des courbe elliptiques sur $\Q$) not\'e
$$S=\Sh(\GL_2,\Hb^\pm).$$
Si $x\in S(\C)$ est un point sp\'ecial correspondant \`a une courbe elliptique
\`a multiplication complexe par un corps quadratique totalement imaginaire
$F/\Q$, on lui associe une sous-vari\'et\'e
de Shimura profinie d\'efinie sur $F$
$$\Sh(T,\{x\})\subset S_F$$
avec $T=\Res_{F/\Q}\G_m$.

On a une action naturelle $\nat$ de $\Gal(\qb/F)$ sur $\Sh(T,\{x\})(\qb)$
donn\'ee par la restriction de l'action naturelle sur $S(\qb)$. On a aussi
une action naturelle de $T(\A_f)/T(\Q)$ sur $\Sh(T,\{x\})(\qb)$.
Ce groupe \'etant le groupe des composantes connexes du groupe de classes
d'id\`eles de $F$, on a un morphisme de r\'eciprocit\'e du corps de classe
$$\artin^{-1}:\Gal(\qb/F)\to T(\A_f)/T(\Q)$$
et donc une nouvelle action de Galois sur $\Sh(T,\{x\})(\qb)$
not\'ee $\rec$ (donn\'ee par l'inverse de l'action $\artin^{-1}$).

Le th\'eor\`eme principal de la multiplication complexe implique que ces deux
actions sont les m\^emes, i.e.
$$\nat=\rec.$$

On peut montrer que ce r\'esultat, bien que visuellement loin du probl\`eme
de Hilbert, en donne la r\'eponse pour les corps quadratiques (la fonction
$j$ et les fonctions elliptiques sont les fonctions dont on prend des
valeurs sp\'eciales car elles permettent de calculer les corps de d\'efinition
des classes d'isomorphismes de courbes elliptiques et de leurs points de torsion).

En rang sup\'erieur, la th\'eorie de la multiplication complexe se formule
de la m\^eme mani\`ere en regardant les points sp\'eciaux dans les espaces
de modules de vari\'et\'es ab\'eliennes principalement polaris\'ees
$$\Sh(\GSp_{2n},\Sc^\pm),$$
qui sont aussi munis de deux actions, l'une naturelle et l'autre donn\'ee
par le corps de classe. Le th\'eor\`eme principal de la multiplication complexe
donne l'\'egalit\'e de ces deux actions.

Les vari\'et\'es de Shimura (connexes) finies sont des (limites projectives de)
quotients du type
$$\Gamma\backslash G(\R)/K$$
avec $G/\Q$ groupe r\'eductif connexe, $\Gamma\subset G(\Q)$ sous-groupe arithm\'etique
et $K\subset G(\R)$ sous-groupe compact maximal dans $G^{ad}(\R)$.

Pour les compactifier, on rajoute des composantes de bords rationnelles
correspondant \`a des sous-groupes paraboliques rationnels $P\subset G$.

On peut aussi compactifier l'espace sym\'etrique $G(\R)/K$ en rajoutant
des composantes du type $G(\R)/P(K)$ o\`u
$P(K)=M(K)AN$ est un parabolique r\'eel (voir \cite{BoJi}). On va s'int\'eresser aux telles
composantes qui viennent d'un parabolique rationnel $P$. On peut voir
les quotients
$$\Gamma\backslash G(\R)/P(K)$$
comme des \'epaississements du bord rationnel, d\'ecrivant certaines d\'eg\'en\'erescences
irrationnelles des structures de Hodge param\'etr\'ees par la vari\'et\'e de Shimura de
d\'epart.
On va montrer (sur des exemples) que les ensembles
$\Gamma\backslash G(\R)/P(K)$ et surtout leurs rev\^etements
$$\Sh^\pm=\Gamma\backslash G(\R)/M(K)A^+$$
et leur interpr\'etation modulaire en termes d'alg\`ebre lin\'eaire
ont un int\'er\^et pour la g\'en\'eralisation du r\^eve de Manin en rang sup\'erieur.

Plus pr\'ecis\'ement, on d\'efinira des ensembles de Shimura
$$\Sh(\GSp_{2n},\Rc^\pm)$$
et des points sp\'eciaux dans ces ensembles et on munira les sous-ensembles
sp\'eciaux correspondants d'une action de Galois $\rec$ donn\'ee
par la th\'eorie du corps de classe.

Le r\^eve de Manin en rang sup\'erieur peut alors (en premi\`ere approche na\"ive)
se formuler en disant qu'il
existe une bonne notion de vari\'et\'e ab\'elienne non commutative sur un
corps de nombres, \`a multiplication par un corps quadratique $F/E$ 
(avec $E/\Q$ un corps totalement r\'eel)
telle que
\begin{enumerate}
\item les p\'eriodes (K-theorie, cohomologie cyclique,\dots)
soient les objets d'alg\`ebre lin\'eaire d\'ecrits
dans ce document,
\item les modules alg\'ebriques soient fortement li\'es aux nombres de
Stark,
\item les sous-espaces sp\'eciaux soient naturellement munis d'une action de Galois $\nat$,
\end{enumerate}
et qu'on ait en plus
$$\rec=\nat.$$
On peut appeler ce r\^eve en rang sup\'erieur le r\^eve de multiplication
quadratique car les corps quadratiques y appara\^issant ne sont pas
n\'ecessairement totalement imaginaires ni totalement r\'eels.

Bien que ce r\^eve soit pure sp\'eculation, il pose quelques questions pr\'eliminaires,
qui ont un int\'er\^et ind\'ependant et dont voici quelques exemples:
\begin{itemize}
\item Est-il possible de donner
une d\'efinition ad\'elique des nombres de Stark pour les corps quadratiques r\'eels?
\item Est-ce que les (g\'en\'erateurs des) sous-cat\'egories stables par extensions des cat\'egories
d\'eriv\'ees de faisceaux coh\'erents
sur une courbe elliptique d\'efinie sur $\Q$
(obtenues en g\'en\'eralisant l\'eg\`erement la construction
de Polishchuk/Schwarz \cite{Poli2}) donn\'ees par nos analogues r\'eels de structures
complexes sont d\'efinies sur le corps de classe du corps de multiplication r\'eelle?
\end{itemize}

%************************************************************************
\section{Remerciements}
Ce document contient essentiellement les notes d'un expos\'e fait \`a Goettingen pour
la conf\'erence ``G\'eom\'etrie diophantienne'' organis\'ee par Yuri Tschinkel. Cet
expos\'e a aussi \'et\'e fait dans plusieurs autres endroit d'avril \`a juin 2004:
les s\'eminaires de g\'eom\'etrie alg\'ebrique de Rennes, K-th\'eorie et g\'eom\'etrie non
commutative \`a l'IHP, math\'ematiques pures de l'ENS Lyon, g\'eom\'etrie non commutative et
th\'eorie des nombres II au Max Planck Institut f\"ur Mathematik de Bonn.
Je remercie mes h\^otes math\'ematiques
(Pierre Berthelot, Max Karoubi, Etienne Ghys, Matilde Marcolli et Yuri Tschinkel) de
m'avoir offert l'opportunit\'e d'exposer ces travaux, qui consistent essentiellement en un effort
de compr\'ehension de l'aspect alg\`ebre lin\'eaire du probl\`eme de multiplication r\'eelle de Manin
en rang sup\'erieur du point de vue des groupes alg\'ebriques. J'esp\`ere que ce point
de vue sera utile aux g\'eom\`etres non commutatifs d\'esireux d'aborder cette question.
Ce travail porte l'empreinte des discussions \'eclairantes que j'ai eues avec Matilde
Marcolli, et pour lesquelles je la remercie. Je remercie aussi Gabor Wiese pour son aide sur la
multiplication complexe. Ce travail \`a \'et\'e fait sur les fonds du r\'eseau
RTN ``K-theory and algebraic groups'', \`a Regensburg. Il suit et compl\^ete les r\'esultats
de l'article \cite{Famille-univ-propre}, mais peut \^etre lu ind\'ependamment de ce dernier.

%************************************************************************
\section{Rappels: corps de classe et vari\'et\'es de Shimura}
%************************************************************************
\subsection{Corps de classe}
On note $\hat{\Z}=\prod\Z_p$, $\A_f\cong \hat{\Z}\otimes_\Z\Q$ et
$\A=\R\times\A_f$
\footnote{On peut aussi d\'ecrire les ad\`eles comme la r\'eunion
des $\A_S:=\R\times\prod_{p\in S}\Q_p\times\times \prod_{p\notin S}\Z_p$
avec $S$ fini. Les ouverts de $\A_S$ sont les
$U\times\prod_{p\notin S}\Z_p$
avec $U$ ouvert de $\R\times\prod_{p\in S}\Q_p$.}.
Si $F/\Q$ est un corps de nombres, on note
$\A_{f,F}=\A_f\otimes_\Q F$ et $\A_F=\A\otimes_\Q F$.

Si $C$ est un groupe topologique commutatif, on note
$$\pi_0(C):=\limproj C/K$$
o\`u la limite projective est prise sur les sous-groupes
$K\subset C$ ferm\'es d'indice fini
\footnote{Cette d\'efinition un peu tordue s'identifie dans le cas qu'on
regarde d'apr\`es Deligne
au $\pi_0$ muni de la topologie quotient. En effet, ce groupe de composantes
connexes des classes d'id\`eles est profini
car compact et totalement discontinu. Il est compact car quotient d'un groupe compact:
$\pi_0(\Ib^1/\Q^\times)$}.

Soit $C_F:=\A_F^\times/F^\times$
\footnote{La topologie induite par la topologie ad\'elique sur les points
d'une vari\'et\'e affine est bonne si cette vari\'et\'e est ferm\'ee.
On prend donc la topologie induite par $\A^\times\mapsto \Ab\times\Ab,
x\mapsto (x,x^{-1})$ car $\G_{m,\Q}\to \A^1_\Q$ (espace affine) n'est
pas ferm\'ee.}.
La th\'eorie du corps de classe nous dit qu'il existe un morphisme
naturel
$$\artin:C_F\to\Gal(F^{ab}/F),$$
qui induit un isomorphisme
$\artin:\pi_0(C_F)\to\pi_0(\Gal(F^{ab}/F))=\Gal(F^{ab}/F)$
et que ce morphisme permet de d\'ecrire toutes les extensions ab\'eliennes
de $F$. On note $\rec:\Gal(F^{ab}/F)\to\pi_0(C_F)$ son inverse.

Soit $T=\Res_{F/\Q}\G_m$ le groupe $F^\times$ vu comme groupe alg\'ebrique
sur $\Q$. Alors on a $C_F=T(\Q)\backslash T(\A)$
et $\pi_0(T(\Q)\backslash T(\A))$ est un
groupe profini qui peut s'\'ecrire comme
$$\limproj_K T(\Q)\backslash T(\A)/T(\R)^+\times K$$
o\`u $K$ parcourt les sous-groupes compacts ouverts de $T(\A_f)$
(Voir Deligne, Vari\'et\'es de Shimura, \cite{De4}).

Un th\'eor\`eme de Chevalet nous dit que
$$\limproj T(\Q)\backslash T(\A_f)/K\cong \overline{T(\Q)}\backslash T(\A_f).$$

Supposons que $F$ est une extension quadratique de $\Q$.
Si $F$ est un corps totalement imaginaire, $T(\R)^+=T(\R)=\C^\times$
et de plus $T(\Q)$ est ferm\'e dans $T(\A_f)$
d'o\`u on d\'eduit
$$\pi_0(T(\Q)\backslash T(\A))=T(\Q)\backslash T(\A_f).$$

Si $F/\Q$ est une extension quadratique r\'eelle,
la composante $T(\R)/T(\R)^+$ est isomorphe \`a $(\Z/2\Z)^2$ et il
ne faut pas l'oublier.

L'espace topologique naturel dans lequel on peu plonger $T(\R)/T(\R)^+$
est l'espace
$$\Rc^\pm:=\GL_2(\R)/T(\R)^+$$
des g\'eod\'esiques orient\'ees sur le demi-plan de poincar\'e.

%************************************************************************
\subsection{Vari\'et\'e de Shimura}
Une \emph{pr\'e-donn\'ee de Shimura} est un couple $(G,X)$ avec $G$ un groupe
alg\'ebrique r\'eductif sur $\Q$ et $X$ un $G(\R)$-espace topologique (ou lisse)
\`a gauche.

Si $(G,X)$ est une pr\'e-donn\'ee de Shimura et $K\subset G(\A_f)$ est
un sous-groupe compact ouvert, on peut lui associer
l'\emph{espace de Shimura fini}
$$\Sh_K(G,X):=G(\Q)\backslash X\times G(\A_f)/K.$$
On appelle \emph{espace de Shimura} la limite projective ensembliste
$$\Sh(G,X):=\limproj_K\Sh_K(G,X).$$

Soit $\Sb:=\Res_{\C/\R}\G_m$.
Une donn\'ee de Shimura est un couple $(G,X)$ avec $G$ r\'eductif connexe
sur $\Q$ et $X\subset\Hom(\Sb,G_\R)$ une $G(\R)$-classe de conjugaison
v\'erifiant quelques axiomes fondamentaux suppl\'ementaires (voir Deligne \cite{De4})
qui impliquent notamment que
$\Sh(G,X)$ est une vari\'et\'e alg\'ebrique quasi-projective (th\'eor\`eme de
Bailly-Borel) qui admet
un mod\`ele canonique sur un corps de nombres (Travaux de Shimura, Deligne, Milne, Shi).
Ces deux th\'eor\`emes a eux seuls montrent toute la puissance des axiomes de base
des vari\'et\'es de Shimura.

Par exemple, si $h_0:\Sb\to\GL_{2,\R}$ est le morphisme donn\'e sur les
points r\'eels par
$z=a+ib\mapsto \left(\begin{smallmatrix}
a & b\\
-b & a\end{smallmatrix}\right)$ et $\Hb^\pm$ est la classe de
$\GL_2(\R)$-conjugaison de $h_0$ alors $(\GL_2,\Hb^\pm)$ est
la donn\'ee de Shimura modulaire dont on a parl\'e dans l'introduction.

Nous allons maintenant nous int\'eresser \`a des pr\'e-donn\'ees de Shimura qui
ne v\'erifient pas les axiomes de donn\'ees de Shimura.

%************************************************************************
\section{G\'eod\'esiques orient\'ees et corps quadratiques r\'eels}
Rappelons par un tableau une analogie entre la courbe modulaire et l'espace des
g\'eod\'esiques sur ic\`ele, essentiellement due \`a Gauss.

\begin{defin}
Soit $M$ un $\Z$-module libre de rang $2$, pour lequel on ne choisit pas
\`a priori de base. Une \emph{structure complexe} sur $M$ est une d\'ecomposition
$M_\C=F\oplus \bar{F}$ de $M_\C$ en deux sous-espaces complexes conjugu\'es.
Un \emph{lilas} sur $M$ est la donn\'ee d'une d\'ecomposition $M_\R=F\oplus \tilde{F}$
de $M_\R$ en somme directe de deux droites et d'orientations sur ces deux droites
donn\'ees par le choix de deux demi-droites $F^+$ et $\tilde{F}^+$ sur ces derni\`eres.
\end{defin}

L'analogie qui nous int\'eresse se fait entre les espaces de modules de ces deux types
d'objets d'alg\`ebre lin\'eaire. Elle comporte trois volets: \emph{alg\`ebre lin\'eaire},
\emph{g\'eom\'etrie analytique} et \emph{g\'eom\'etrie alg\'ebrique}.
$$
\begin{array}{c}
\textbf{Alg\`ebre lin\'eaire}\\
\begin{array}{|c|c|c|}\hline
 & \textrm{structures complexes} & \textrm{lilas}\\ \hline
\textrm{en g\'en\'eral} &
M_\C=F\oplus \bar{F}&
M_\R=F\oplus \tilde{F}\textrm{ et }F^+\subset F,\tilde{F}^+\subset \tilde{F}\\ \hline
\textrm{un exemple sp\'ecial} &
M=\Z[\sqrt{-2}] & M=\Z[\sqrt{2}]\\ \hline
\begin{array}{c}
\textrm{isomorphismes}\\
\textrm{d'alg\`ebres}
\end{array} &
\Z[\sqrt{-2}]\otimes_\Z\C\underset{alg}{\cong}
\C_{\iota_1}\times\C_{\iota_2} &
\begin{array}{c}
\Z[\sqrt{2}]\otimes_\Z\R\underset{alg}{\cong}\R_{\iota_1}\times \R_{\iota_2},\\
\R^+_{\iota_1},\R^+_{\iota_2}
\end{array}\\ \hline
\end{array}
\end{array}
$$
Dans la derni\`ere ligne du tableau, $\iota_1$ et $\iota_2$ sont les deux plongements
de l'alg\`ebre consid\'er\'ee,
complexes pour la colonne de gauche et r\'eels pour la colonne de droite.
L'orientation sur les droites r\'eelles est obtenue en prenant la composante
connexe de l'identit\'e dans le groupe des inversibles $\R^*$.
\begin{defin}
Un \emph{morphisme} entre deux structures complexes (resp. lilas) est la donn\'ee d'un
morphisme $f:M_1\to M_2$ entre les r\'eseaux sous-jacents tel que
$f_\C(F_1)\subset F_2$ (resp. $f_\R(F_1^+)\subset F_2^+$ et
$f_\R(\tilde{F}_1^+)\subset \tilde{F}_2^+$).
\end{defin}
\begin{defin}
Une \emph{structure de niveau $N$} sur une structure complexe (resp. un lilas) de r\'eseau
sous-jacent $M$ est un isomorphisme $M\otimes_\Z\Z/N\Z\overset{i}{\to} (\Z/N\Z)^2$.
\end{defin}
Nous nous int\'eressons maintenant \`a la g\'eom\'etrie des espaces de modules de ces
objects d'alg\`ebre lin\'eaire.
$$
\begin{array}{c}
\textbf{G\'eom\'etrie analytique}\\
\begin{array}{|c|c|c|}\hline
 & \textrm{classique} & \textrm{dynamique}\\ \hline
\begin{array}{c}
\textrm{espace}\\
\textrm{de modules}
\end{array} &
\GL_2(\Z)\backslash\GL_2(\R)/\C^\times&
\GL_2(\Z)\backslash\GL_2(\R)/(\R^{+\times})^2\\ \hline
\textrm{g\'eom\'etriquement} &
X=\PGL_2(\Z)\backslash\Hb&
\{\textrm{g\'eod\'esiques sur }$X$\}\\ \hline
\begin{array}{c}
\textrm{tour modulaire}\\
\textrm{(tous niveaux)}
\end{array}&
\limproj \Gamma(N)\backslash\GL_2(\R)/\C^\times &
\limproj \Gamma(N)\backslash\GL_2(\R)/(\R^{+\times})^2
\end{array}\\ \hline
\end{array}
$$
Les groupes $\Gamma(N)$ consid\'er\'es dans la limite projective sont les
usuels groupes de congruence de niveau $N$, donn\'es par la suite exacte
$$1\to \Gamma(N)\to\GL_2(\Z)\to \GL_2(\Z/N\Z)\to 1.$$
Cet aspect g\'eom\'etrique est celui qui semble bien connu par Gauss.
Une des mani\`eres d'aborder la question de multiplication r\'eelle de Manin
est de chercher \`a remplir le tableau \`a trous suivant.
$$
\begin{array}{c}
\textbf{G\'eom\'etrie alg\'ebrique}\\
\begin{array}{|c|c|c|}\hline
 & \textrm{classique} & \textrm{quantique}\\ \hline
\textrm{structures lin\'eaires} &
M_\C=F\oplus \bar{F}&
M_\R=F\oplus \tilde{F}\textrm{ et }F^+\subset F,\tilde{F}^+\subset \tilde{F}\\ \hline
\textrm{r\'ealisation g\'eom\'etrique} &
E=M\backslash M_\C/\bar{F}\subset \Pb^2(\C) &
?\\ \hline
\textrm{par exemple} &
E=\C/\Z[\sqrt{-2}]\subset \Pb^2(\C) &
? \sim \R/\Z[\sqrt{2}]\\ \hline
\end{array}
\end{array}
$$

On va maintenant expliquer plus en d\'etails la deuxi\`eme colonne de nos tableaux
\emph{alg\`ebre lin\'eaire} et \emph{g\'eom\'etrie analytique} en gardant cette analogie en t\^ete.

Soit $h_1:\G_{m,\R}^2\to \GL_{2,\R}$ le morphisme donn\'e par
$(x,y)\mapsto\diag(x,y)$ et $\Rc\subset\Hom(\G_{m,\R}^2,\GL_{2,\R})$
sa $\GL_2(\R)$-classe de conjugaison.
Le stabilisateur de $h_1$ est le tore diagonal $T(\R)$ de $\GL_2(\R)$ donc
$$\Rc\cong\GL_2(\R)/T(\R).$$

Les points de $\Rc$ correspondent aux couples de deux droites propres
pour $x$ et $y$. On peut les dessiner comme des points sur le bord
$\Pb^1(\R)$ du disque de poincar\'e. On munit maintenant ces points d'un sens
(arriv\'ee ou d\'epart).

Soit $\Rc^\pm$ l'espace des couples $(h,c)$ avec $h\in \Rc$
et $c\in\pi_0(\im(h(\R)))$. On a alors
$$\Rc^\pm\cong (\Z/2\Z)^2\times \Rc$$
et l'action de $\GL_2(\R)$ par conjugaison agit seulement sur le deuxi\`eme
facteur.

\begin{defin}
Le couple $(\GL_2,\Rc^\pm)$ est appel\'e la \emph{donn\'ee de rivage modulaire}.
\end{defin}

\begin{defin}
Si $h\in \Rc$, on appelle \emph{mauvais
\footnote{Les bons groupes de Mumford-Tate sont donn\'es par les enveloppes
r\'eductives et ont un int\'er\^et ind\'ependant pour des questions de
dynamique.}
groupe de Mumford-Tate de $h$} le
plus petit sous-$\Q$-groupe alg\'ebrique $\BMT(h)\subset\GL_2$
contenant sur $\R$ l'image de $h$. Un point $h\in \Rc$ est dit
\emph{sp\'ecial} si $\BMT(h)$ est un tore non scind\'e.
\end{defin}

On montre facilement que si $h$ est sp\'ecial, $\BMT(h)$ peut s'\'ecrire
$\Res_{E/\Q}\G_m$ avec $E/\Q$ quadratique totalement r\'eel
\footnote{Par exemple: consid\'erons $d>1$ un entier positif sans facteur carr\'e.
Soit
$g:=\left(\begin{smallmatrix}
1 & 1\\
\sqrt{d} & -\sqrt{d}
\end{smallmatrix}
\right)$.
En conjuguant $h_0$ par $g$, on obtient la matrice
$h'=\left(\begin{smallmatrix}
a & b\\
db & a
\end{smallmatrix}
\right)$ avec $a=\frac{x+y}{2}$ et $b=\frac{x-y}{2\sqrt{d}}$.
Le groupe des matrices de la forme
$\left(\begin{smallmatrix}
a & b\\
db & a
\end{smallmatrix}
\right)$
avec $a$ et $b$ rationnels est un groupe
alg\'ebrique sur $\Q$ qui est le groupe de Mumford-Tate de $h'$.
C'est en fait $\Res_{\Q(\sqrt{d})/\Q}\G_m$.}.

On peut ainsi compl\'eter notre analogie par un tableau intitul\'e
\emph{structure sp\'eciales} qui donne les id\'ees clefs de la question
de multiplication r\'eelle.
$$
\begin{array}{c}
\textbf{Structures sp\'eciales}\\
\begin{array}{|c|c|c|}\hline
 & \textrm{multiplication complexe} & \textrm{multiplication r\'eelle}\\ \hline
\textrm{structures lin\'eaires} &
\begin{array}{c}
M_\C=F\oplus \bar{F}\\
\End(M,F,\bar{F})\otimes_\Z\Q\\
\textrm{quadratique imaginaire}
\end{array} &
\begin{array}{c}
M_\R=F\oplus \tilde{F}\textrm{ et }F^+\subset F,\tilde{F}^+\subset F\\
\End(M,F,\tilde{F})\otimes_\Z\Q\\
\textrm{quadratique r\'eel}
\end{array}\\ \hline
\textrm{par exemple} &
M=\Z[\sqrt{-2}] &
M=\Z[\sqrt{2}]\\ \hline
\textrm{g\'eom\'etrie alg\'ebrique} &
\begin{array}{c}
\textrm{courbes elliptiques}\\
\textrm{\`a multiplication complexe}
\end{array} &
?\\ \hline
\textrm{structures de niveau} &
\textrm{points de torsion}&
?\\ \hline
\textrm{modules alg\'ebriques} &
\begin{array}{c}
\textrm{valeurs de fonction $j$}\\
\textrm{et fonctions elliptiques}
\end{array} &
?\sim\textrm{nombres de stark}\\ \hline
\end{array}
\end{array}
$$

On peut dessiner sur le disque de poincar\'e
la g\'eod\'esique liant $\sqrt{2}$ et $-\sqrt{2}$
qui est une g\'eod\'esique sp\'eciale.

Si $h\in \Rc$ est un tel point sp\'ecial de groupe de Mumford-Tate
$T=\Res_{E/\Q}\G_m$,
on lui associe la sous-donn\'ee
$$s:(T,Z=\pi_0(T(\R)))\to (\GL_2,\Rc^\pm).$$

Comme vu pr\'ec\'edemment, on a
$$\Sh(T,Z)=\pi_0(T(\Q)\backslash T(\A))$$
s'identifie au groupe des composantes connexes du groupe de classes d'id\`eles.

On munit donc l'image de $\Sh(s)$ d'une action not\'ee
$\rec$ de $\Gal(E^{ab}/E)$ donn\'ee par l'inverse de l'isomorphisme
du corps de classe.
\begin{reve}[``Alterstraum'' de Manin: multiplication r\'eelle]
On esp\`ere que cette action
peut-\^etre d\'ecrite de mani\`ere plus naturelle gr\^ace \`a
une interpr\'etation modulaire de $\Sh(T,Z)$ en termes de g\'eom\'etrie
alg\'ebrique non commutative, li\'ee aux nombres de Stark pour ce corps quadratique.
\end{reve}

%************************************************************************
\section{Dans les espaces de modules de vari\'et\'es ab\'eliennes}
On souhaite comprendre ce que le r\^eve de Manin signifie pour un corps comme
$\Q(\sqrt[4]{2})/\Q(\sqrt{2})$. On remarque qu'on a
$$\Q(\sqrt[4]{2})\otimes_\Q\R
\underset{\R-alg}{\cong} \C\times \R\times\R$$
et ce corps quadratique n'est ni totalement r\'eel, ni totalement imaginaire.

Soit
$\Db_k\subset(\G_{m,\R}^2)^k$ le sous-groupe des $(x_i,y_i)_{i=1,\dots,k}$
tels que $x_iy_i=x_jy_j$ pour tout $i,j\in \{1,\dots,k\}$.
On remarque que $\Db_k$ peut \^etre vu comme le tore maximal
de $\GSp_{2k,\R}$.
Soit $\Tb_k\subset \Sb\times \Db_k$ le sous-groupe des $(z,(x_i,y_i))$
tels que
$z\bar{z}=x_iy_i$ pour tout $i\in\{1,\dots,k\}$.
Soit $\Tb:=\limproj \Tb_k$ o\`u les projections $\Tb_{k+1}\to\Tb_k$
sont donn\'ees par oubli du dernier facteur $\G_{m,\R}^2$.
On note $\pi_k:\Tb\to\Tb_k$ la projection naturelle.

On a un morphisme diagonal $w:\G_m\to\Tb$ dit de poids et un autre morphisme
$\mu:\G_{m,\C}\to\Tb_\C$ dit de Hodge donn\'e par
$x\mapsto ((x,1),\dots,(x,1))$.

Soit $G_n\subset\GL_2^n$ le sous-groupe des matrices $(g_i)$ telles que
$\det(g_i)=\det(g_j)$ pour tout $i,j$.

Soit
$$f_n:G_n\to\GSp_{2n}$$
le morphisme donn\'e pour
$(g_i=\left(\begin{smallmatrix}
a_i & b_i\\
c_i & d_i\end{smallmatrix}\right))$ par la matrice
$$
\left(\begin{array}{cc}
\diag(a_i) & \diag(b_i)\\
\diag(c_i) & \diag(d_i)
\end{array}\right)
$$

On note
$h_0:\Sb\to\GL_{2,\R}$ le morphisme standard donn\'e sur les
points r\'eels par
$z=a+ib\mapsto \left(\begin{smallmatrix}
a & b\\
-b & a\end{smallmatrix}\right)$
et
$h_1:\G_{m,\R}^2\to\GL_{2,\R}$ le morphisme diagonal
donn\'e par $(x,y)\mapsto \diag(x,y)$.

Soit $k_0+k_1=n$ une partition de $n$.
On note $h_{k_0,k_1}:\Tb_{k_1}\to G_n\subset\GL_{2,\R}^n$ le morphisme
donn\'e par
$$
h_{k_0,k_1}:=
(\underset{k_0}{\underbrace{h_0,\dots,h_0}},
\underset{k_1}{\underbrace{h_1,\dots,h_1}}).$$

Soit $h_{S,k_0,k_1}:\Tb\to \GSp_{2n}$ le morphisme donn\'e par
$h_{S,k_0,k_1}=f_n\circ h_{k_0,k_1}\circ\pi_{k_1}$ et $\Rc_{k_0,k_1}$
la $\GSp_{2n}(\R)$-classe de conjugaison de ce morphisme dans
$\Hom(\Tb,\GSp_{2n})$. On remarque que pour $h\in \Rc_{k_0,k_1}$, le morphisme
$h\circ w:\G_{m,\R}\to \GSp_{2n}$ est simplement le plongement diagonal et est
donc ind\'ependant de $h$.

Dans le cas $k_1=0$, on retrouve l'espace de Siegel classique $\Sc^\pm$
et dans le cas $k_0=0$, on trouve
$$\Rc_{0,n}\cong \GSp_{2n}(\R)/\Db_n(\R).$$

Rappelons que le couple $(\GSp_{2n},\Sc^\pm)$ est appel\'ee \emph{donn\'ee de Shimura
de Siegel}. La vari\'et\'e de Shimura correspondante est l'espace de module
des vari\'et\'es ab\'eliennes \`a structure de niveau infinie.

On note $\Rc_{k_0,k_1}^\pm$ l'espace des morphismes $h$ munis d'une composante
connexe $c\in\pi_0(\im(h(\Tb(\R))))$.

\begin{defin}
Le couple $(\GSp_{2n},\Rc_{k_0,k_1}^\pm)$ est appel\'e \emph{donn\'ee de rivage
de type $(k_0,k_1)$} de la donn\'ee de Siegel $(\GSp_{2n},\Sc^\pm)$.
\end{defin}

\begin{defin}
Si $h\in \Rc_{k_0,k_1}$, on appelle \emph{mauvais
\footnote{Les bons groupes de Mumford-Tate sont donn\'es par les enveloppes
r\'eductives. Il semble plausible que ces groupes soient d\'efinis de mani\`ere
unique (c'est le cas si $k_0=0$ ou $k_1=0$). Il serait int\'eressant de conna\^itre
leur signification pour des questions de dynamique.}
groupe de Mumford-Tate de $h$} le
plus petit sous-$\Q$-groupe alg\'ebrique $\BMT(h)\subset\GSp_{2n}$
contenant sur $\R$ l'image de $h$.
Un point $h\in \Rc_{k_0,k_1}$ est dit \emph{sp\'ecial} si son groupe de
Mumford-Tate est un tore $T$ tel que $T/w(\G_m)$ soit anisotrope sur $\Q$
\footnote{Cette condition est automatiquement v\'erifi\'ee sur $\R$ dans
le cas des vari\'et\'es ab\'eliennes car $h(i)$ est une involution de Cartan
de $(T/w(\G_m))_\R$.}.
\end{defin}

Si $h$ est un point sp\'ecial, on a
$$T(\R)/T(\R)^+=h(\Tb(\R))/h(\Tb(\R))^+$$
d'o\`u un plongement
$$T(\R)/T(\R)^+\subset \Rc^\pm_{k_0,k_1}$$
donn\'e par $c\mapsto (h,c)$.

On peut maintenant refaire la construction de Deligne dans ce cadre.
On a un morphisme
$\mu_h:=h_\C\circ\mu:\G_{m,\C}\to T_\C$
dit de Hodge et le corps reflex $E(h)$ de $h$ est le corps de d\'efinition
de ce morphisme. Si $F=E(h)$, on a
$$\mu_h:\G_{m,F}\to T_F$$
et on construit
$$\NR(\mu_h):\Res_{F/\Q}\G_m\to\Res_{F/\Q}T_F\overset{\Nm}{\to} T.$$

On a
$$\pi_0(\NR(\mu_h)):\pi_0(C_F)\to \pi_0(T(\Q)\backslash T(\A)).$$

Soit $h$ un point sp\'ecial de groupe de Mumford-Tate $T$,
on munit canoniquement le sous-espace
$$\Sh(T,Z)\subset\Sh(\GSp_{2n},\Rc_{k_0,k_1}^\pm)$$
avec $Z=\pi_0(T(\R))$
d'une action de Galois not\'ee $\rec$ donn\'ee par
$$g:x\mapsto (\pi_0(\NR(\mu_h))\circ \rec(g))\cdot x.$$

\begin{reve}[de multiplication quadratique]
On esp\`ere que cette action de Galois peut-\^etre construite de mani\`ere
plus naturelle en interpr\'etant $\Sh(T,Z)$ comme un espace de modules
d\'efini sur un corps de nombres d'analogues non commutatifs des vari\'et\'es
ab\'eliennes dont les modules alg\'ebriques seraient li\'es aux nombres de Stark
pour $F$.
\end{reve}

Regardons notre exemple de d\'epart de ce point de vue.
Notons $F'=\Q(\sqrt[4]{2})$ et $E=\Q(\sqrt{2})$. Le choix d'une $E$-base de $F'$
nous donne un plongement $\Res_{F'/\Q}\G_m\subset \Res_{E/\Q}\GL_2$.
On note $T\subset G$ l'inclusion correspondante des sous-groupes donn\'es par
les matrices de d\'eterminant dans $\G_{m,\Q}$. On remarque que
$T\cong \Res_{F'/\Q}\G_m/(\Res_{E/\Q}\G_m)^{(1)}$.
On peut plonger $G$ dans $\GSp_{4,\Q}$.
On a $T_\R\cong \Tb_1$ et on note $h:\Tb\to \GSp_{4,\R}$ le morphisme
correspondant. Alors $h\in \Rc_{1,1}$ et $T$ est clairement le groupe de Mumford-Tate de $h$.
Le morphisme $\mu_h:\G_{m,\C}\to T_\C$ est d\'efini sur le corps $F=\Q(i\sqrt[4]{2},\sqrt[4]{2})$.

%************************************************************************
\section{Litt\'erature}
La litt\'erature regarde essentiellement des cas de rang $1$.
Les gens s'int\'eressent au bord irrationnel
$$\PGL_2(\Z)\backslash\Pb^1(\R)$$
qui a l'air d'\^etre un espace moins naturel que le rivage
$$\GL_2(\Z)\backslash\GL_2(\R)/{\R^{+\times}}^2$$
de notre point de vue de la th\'eorie du corps de classe et
des groupes alg\'ebriques.

On d\'ecoupe la litt\'erature en deux approches:
\begin{itemize}
\item l'approche \emph{p\'eriodes} qui part des objets d'alg\`ebre
lin\'eaire et construit des objets de g\'eom\'etrie non commutative
correspondants,
\item l'approche \emph{nombres de Stark} qui part des nombres de Stark
et construit des objets de g\'eom\'etrie non commutative correspondants.
\end{itemize}

Le r\^eve est qu'il existe un bon point de rencontre entre ces deux approches.

La litt\'erature s'\'enum\`ere ainsi.
\begin{enumerate}
\item Connes, Manin, Marcolli \cite{Connes3,Manin3,Manin-Marcolli2}:
approche \emph{periodes}: les aspects
analytiques peuvent \^etre donn\'es par les tores non commutatifs
($K$-theorie, $\HC$).
\item Polishchuk \cite{Poli1}: approche \emph{periodes}: les tores non commutatifs analytiques
sont des ``vari\'et\'es projectives''.
Annonce en rang sup\'erieur.
\item Manin \cite{Manin3}: approche \emph{nombres de Stark}: fonctions theta non commutatives et
nombres de Stark.
\item Darmon: autre approche: points de Hegner-Stark sur les courbes elliptiques sur $\Q$
conjecturalement d\'efinis sur des corps de classe.
\item Connes-Marcolli-Ramachandran (non publi\'e): \'etats KMS et multiplication complexe. 
\end{enumerate}

\appendix
\section{Rappels sur la multiplication complexe}
Ce paragraphe est essentiellement le r\'esultat de mes discussions avec Gabor Wiese
que je remercie pour son aide.
Soit $F$ un corps de nombres.
Soit $\Ell(F)_\Q$ la cat\'egorie dont les objets sont les courbes elliptiques
sur $F$ et les morphismes sont donn\'es par $\Hom(E_1,E_2)\otimes_\Z\Q$.
On note une courbe elliptique $E$ vue comme objet de $\Ell(F)_\Q$ comme
$E\otimes\Q$.

Pour $E/F$ une courbe elliptique, on note
$$T(E):=\limproj_n E[n](\bar{F})\textrm{ et }V(E):=T(E)\otimes_\Z\Q$$
avec $E[n]$ le noyau de la multiplication par $n$ dans $E$.

On note $\A_f:=\hat{\Z}\otimes_\Z\Q$ avec $\hat{\Z}$ le compl\'et\'e
de $\Z$ pour les sous-groupes d'indice fini.
On sait que $T(E)\cong\hat{\Z}^2$ et donc
$$V(E)\cong\A_f^2.$$

Soit $F/\Q$ un corps quadratique totalement imaginaire. Une courbe elliptique
est dite \`a multiplication complexe par $F$ si on se fixe un morphisme
injectif $F\to\End(E\otimes\Q)$. Si on note $\A_{f,F}:=\A_f\otimes_\Q F$,
on peut montrer que $V(E)$ est un $\A_{f,F}$-module de rang $1$.

On note $S_F(\qb)$ l'espace des classes d'isomorphismes de triplets
$$(E\otimes\Q,\alpha:F\to\End(E\otimes\Q),\psi:\A_{f,F}\to V(E))$$
avec $\psi$ un isomorphisme de $\A_{f,F}$-modules et $\alpha$ une injection.

On a une action naturelle de $\A_{f,F}^\times$ sur $S_F(\qb)$
et on peut montrer que le stabilisateur d'un points est $F^\times$,
ce qui fait de
$S_F(\qb)$ un $C_F:=\A_{f,F}^\times/F^\times$-torseur.
D'autre part, $S_F(\qb)$
est muni d'une \underline{action naturelle}
de $\Gc_F^{ab}:=\Gal(F^{ab}/F)$ et ces deux actions commutent,
i.e. $(\sigma\cdot x)\times y=\sigma\cdot(x\times y)$.

Fixons un point $x\in S_F(\qb)$, donc un isomorphisme
$C_F\cong S_F(\qb)$.
Ceci nous permet de d\'efinir un morphisme
\footnote{On a $\Phi(\sigma)\times\Phi(\tau)=
(\sigma\cdot 1)\times(\tau\cdot 1)=\sigma\cdot(1\times(\tau\cdot 1)))=
\sigma\cdot(\tau\cdot 1)=(\sigma\tau)\cdot 1$.}
$$\Phi:\Gc_F^{ab}\to C_F$$
par $\Phi(\sigma)=\sigma\cdot 1$.

On a d'autre part l'application d'artin donn\'ee par le corps de
classe $\artin:C_F\to \Gc_F^{ab}$.

\begin{theo}[Principal: Multiplication complexe]
La composition
$$C_F\overset{\artin}{\to} \Gc_F^{ab} \overset{\Phi}{\to} C_F$$
envoie $x\in C_F$ sur $x^{-1}$.
\end{theo}

% *************************
% Appel de la bibliographie
% *************************
\bibliographystyle{hamsplain}
\bibliography{$HOME/travail/fred.bib}

\providecommand{\bysame}{\leavevmode\hbox to3em{\hrulefill}\thinspace}
\begin{thebibliography}{1}

\bibitem{BoJi}
Armand Borel and Ji~Lizhen, \emph{Compactifications of locally symmetric
  spaces}, IAS/Michigan, 2001, Prepublication, june 2001.

\bibitem{Connes3}
Alain Connes, \emph{Noncommutative differential geometry}, Inst. Hautes
  \'Etudes Sci. Publ. Math. (1985), no.~62, 257--360.

\bibitem{De4}
Pierre Deligne, \emph{Vari\'et\'es de {S}himura: interpr\'etation modulaire, et
  techniques de construction de mod\`eles canoniques}, Automorphic forms,
  representations and $L$-functions (Proc. Sympos. Pure Math., Oregon State
  Univ., Corvallis, Ore., 1977), Part 2, Amer. Math. Soc., Providence, R.I.,
  1979, pp.~247--289.

\bibitem{Manin3}
Yuri~I. Manin, \emph{{Real Multiplication and noncommutative geometry}},
  \mbox{arXiv:math.AG/0202109}.

\bibitem{Manin-Marcolli2}
Yuri~I. Manin and Matilde Marcolli, \emph{Continued fractions, modular symbols,
  and noncommutative geometry}, Selecta Math. (N.S.) \textbf{8} (2002), no.~3,
  475--521.

\bibitem{Famille-univ-propre}
Fr\'ed\'eric Paugam, \emph{Three examples of noncommutative boundaries of
  shimura varieties}, Disponible sur http://name.math.univ-rennes1.fr/fpaugam,
  soumis \`a la conf\'erence "Noncommutative geometry and number theory", Bonn,
  ao\^ut 2003. (2004), 1--24.

\bibitem{Poli1}
Alexander Polishchuk, \emph{Noncommutative 2-tori with real multiplication as
  noncommutative projective varieties}, arXiv (2002),
  no.~http://fr.arXiv.org/abs/math.AG/0212306.

\bibitem{Poli2}
\bysame, \emph{Classification of holomorphic vector bundles on noncommutative
  two-tori}, arXiv (2003), no.~http://fr.arXiv.org/abs/math.QA/0308136.

\end{thebibliography}

\end{document}